\documentclass{article}
\usepackage{graphicx}
\usepackage{url}
\usepackage{amsmath,amsthm,amsfonts,amssymb}
\usepackage[colorlinks=true,citecolor=black,linkcolor=black,urlcolor=blue]{hyperref}
\usepackage{cleveref}
\usepackage{float}
\usepackage{enumitem}
\usepackage{todonotes}
\usepackage{xcolor}

\newtheorem{theorem}{Theorem}

\usepackage{xcolor}

\renewcommand{\le}{\leqslant}
\renewcommand{\ge}{\geqslant}

\title{Half-flips are 5-avoidable}
\author{Pascal Ochem\thanks{This work is supported by the ANR project CADO (ANR-24-CE48-3758-01).}\\
\small LIRMM, CNRS, Universit\'e de Montpellier, France}

\date{}

\begin{document}

\maketitle

\begin{abstract}
A word contains a \emph{half-flip} if it contains non-empty factors $uv$ and $vu$ where
$|u|=|v|$. Fici reports a non-constructive proof of the existence of an infinite word
over a finite alphabet avoiding half-flips and asks for the size of the smallest alphabet over which half-flips may be avoided.
Currie and Rampersad have proposed a pure morphic word over 8 letters and a morphic word over 5 letters and conjecture that they avoid half-flips.
We present a pure morphic word over 5 letters that avoids half-flips.
We also show that half-flips with $|u|\ge2$ are 3-avoidable and that half-flips with $|u|\ge4$ are 2-avoidable.
\end{abstract}

\section{Introduction}
A word contains a \emph{half-flip} if it contains non-empty factors $uv$ and $vu$ where $|u|=|v|$.
Let $p=|u|=|v|$ be the \emph{period} of a half-flip.
A $k$-half-flip is a half-flip with period at least $k$.
Following a suggestion of Fici and the work of Currie and Rampersad~\cite{CuRa:2025},
we settle the avoidability of $k$-half-flips.
\begin{theorem}\label{thm-main}{\ }
\begin{enumerate}
\item 1-half-flips are 5-avoidable.
\item 2-half-flips are 3-avoidable.
\item 4-half-flips are 2-avoidable.
\end{enumerate}
\end{theorem}
Using backtracking, it is quite straightforward to check that Theorem 1 is best possible in terms of the alphabet size and the period of avoided half-flip.
\section{Proof of Theorem 1.1}

We use the following 95-uniform morphism:

\begin{align*}
m(\texttt{0})&={\color{red}c}\texttt{{\color{red}13012402301302402301}240241341240230130240230124}\\
m(\texttt{1})&=c\texttt{023012402413012402301302402{\color{blue}34124024134124023013}}\\
m(\texttt{2})&=c\texttt{02341240241301240234130241301240230130240230124}\\
m(\texttt{3})&={\color{red}c}\texttt{{\color{red}13012402301302402301}240241341240230130240234124}\\
m(\texttt{4})&=c\texttt{023412402413012402301302402{\color{blue}34124024134124023013}}
\end{align*}
where $c=\texttt{024130124023013024134124023413024134124023013024}$.
Let us show that $m^\omega(\texttt{0})$ avoids half-flips.
Suppose that $m^\omega(\texttt{0})$ contains a half-flip with period $p$ and that $p$ is minimal,
that is, $m^\omega(\texttt{0})$ contains no half-flip with period smaller than $p$.
A computer check shows that $m^\omega(\texttt{0})$ avoids half-flips with period at most 500.
Thus $p>500$ and $u$ contains at least $\tfrac{501-94}{95}\ge4$ full $m$-images of a letter.
Since $c$ only appears in $m^\omega(\texttt{0})$ as the prefix of the $m$-image of a letter,
all the occurrences of $u$ in $m^\omega(\texttt{0})$ start at the same position modulo 95.
Similarly, all the occurrences of $v$ start at the same position modulo 95.
Let $x_u$ and $x_v$ denote these positions. Since $uv$ is a factor, we have 
$$x_v\equiv x_u+p\pmod{95}.$$
Since $vu$ is a factor, we have 
$$x_u\equiv x_v+p\pmod{95}.$$
By summing up, we obtain
$$x_u+x_v\equiv x_v+x_u+2p\pmod{95}.$$
This gives
$$2p\equiv 0\pmod{95}.$$
Since 95 is odd, $p$ is a multiple of 95. Moreover, $x_u=x_v$.
Let us write $u=u'm(U)u''$ and $v=v'm(V)v''$,
where $0\le|u''|=|v''|=x_u\le94$ and $|u'|=|v'|=95-x_u$. Notice also that $|U|=|V|>1$.
So $m^\omega(\texttt{0})$ contains the factors $uv=u'm(U)u''v'm(V)v''$ and $vu=v'm(V)v''u'm(U)u''$.
Now we consider two cases depending on $x_u$.
\subsection*{case $x_u>70$}
In this case, we have $|u''|>70$. Let us show that there exists a unique letter $\mu$
such that $m(U)u''$ is a prefix of $m(U\mu)$.
Notice that, if $a$ and $b$ are distinct letters, then $m(a)$ and $m(b)$ have distinct prefixes
of length 68, except that $m(\texttt{0})$ and $m(\texttt{3})$ have the same prefix of length 68
depicted in red. So our claim of uniqueness holds, except maybe if $m(U)u''$ is a prefix of both
$m(U\texttt{0})$ and $m(U\texttt{3})$. Let us show that the latter situation does not occur.
Notice that every factor $ab$ of $m^\omega(\texttt{0})$ is such that either $b\equiv a+1\pmod{5}$ or $b\equiv a+2\pmod{5}$.
So, if $w$ is non-empty and $wm$ and $wn$ are distinct factors of $m^\omega(\texttt{0})$, then $m$ and $n$ are consecutive letters,
that is, $m\equiv n+1\pmod{5}$ or $n\equiv m+1\pmod{5}$.
Therefore, $U\texttt{0}$ and $U\texttt{3}$ are not both factors of $m^\omega(\texttt{0})$.
This shows that there is only one way to extend $m(U)u''$ into $m(U\mu)$.
Similarly, there is only one way to extend $m(V)v''$ into $m(V\nu)$.
Now, since $m^\omega(\texttt{0})$ contains $uv=u'm(U)u''v'm(V)v''$, $m^\omega(\texttt{0})$
must contain $u'm(U\mu V\nu)$. By taking the pre-image, $m^\omega(\texttt{0})$ also contains $U\mu V\nu$. 
Symmetrically, $m^\omega(\texttt{0})$ contains $vu=v'm(V)v''u'm(U)u''$, $v'm(V\nu U\mu)$, and $V\nu U\mu$.
So $m^\omega(\texttt{0})$ contains both $U\mu V\nu$ and $V\nu U\mu$.
This is a half-flip with period $\tfrac p{95}$, which contradicts the minimality of $p$.

\subsection*{case $x_u\le70$}
The proof of this case is roughly a mirror image of the previous case.
Here we have $|u'|=95-x_u\ge95-70=25$.
If $a$ and $b$ are distinct letters, then $m(a)$ and $m(b)$ have distinct suffixes
of length 20, except that $m(\texttt{1})$ and $m(\texttt{4})$ have the same suffix of length 20
depicted in blue. Since \texttt{1} and \texttt{4} are not consecutive letters,
there is only one way to extend $u'm(U)$ into $m(\mu U)$
and there is only one way to extend $v'm(V)$ into $m(\nu V)$.
Since $m^\omega(\texttt{0})$ contains $uv=u'm(U)u''v'm(V)v''$, $m^\omega(\texttt{0})$
must contain $m(\mu U\nu V)v''$ and $\mu U\nu V$. 
Symmetrically, $m^\omega(\texttt{0})$ also contains $\nu V\mu U$.
Thus $m^\omega(\texttt{0})$ contains a half-flip with period $\tfrac p{95}$, which is a contradiction.

\section{Proof of Theorems 1.2 and 1.3}
To prove Theorem 1.2 and Theorem 1.3, we will show that $f_3(m^\omega(\texttt{0}))$ avoids 2-half-flips
and $f_2(m^\omega(\texttt{0}))$ avoids 4-half-flips, where the 7-uniform morphisms $f_3$ and $f_2$ are given below.

\noindent
\begin{minipage}[b]{0.4\linewidth}
\centering
$$\begin{array}{l}
 f_3(\texttt{0})=\texttt{0001022}\\
 f_3(\texttt{1})=\texttt{0102122}\\
 f_3(\texttt{2})=\texttt{0211122}\\
 f_3(\texttt{3})=\texttt{0210112}\\
 f_3(\texttt{4})=\texttt{0002122}
\end{array}$$
\end{minipage}
\hspace{2mm}
\begin{minipage}[b]{0.4\linewidth}
\centering
$$\begin{array}{l}
 f_2(\texttt{0})=\texttt{0000001}\\
 f_2(\texttt{1})=\texttt{0101001}\\
 f_2(\texttt{2})=\texttt{0010011}\\
 f_2(\texttt{3})=\texttt{0001111}\\
 f_2(\texttt{4})=\texttt{1011011} 
\end{array}$$
\end{minipage}

Suppose that $f_3(m^\omega(\texttt{0}))$ or $f_2(m^\omega(\texttt{0}))$ contains a forbidden half-flip with period $p$.
A computer check shows that they avoid half-flips forbidden with period at most 500.
Now we need that every not-too-small factor always appears at the same position modulo 7.
Every factor $t$ of $f_3(m^\omega(\texttt{0}))$ such that $|t|\ge8$ contains the factor \texttt{20}.
Since \texttt{20} only appears as the factor overlapping two $f_3$-images, all the occurrences of $t$
have the same position modulo 7.
Also, we can check that for every factor $t$ of $f_2(m^\omega(\texttt{0}))$ of length 6,
all the occurrences of $t$ have the same position modulo 7.
With the same argument as in the previous proof, and since the length 7 is odd,
we obtain that $p$ is a multiple of 7.
Now, if $a$ and $b$ are distinct letters, then $f_3(a)$ and $f_3(b)$ have distinct prefixes (resp. suffixes)
of length 4. Similarly, $f_2(a)$ and $f_2(b)$ have distinct prefixes (resp. suffixes) of length 4.
This is sufficient to invoke the argument of the previous proof that either the prefix or the suffix
extends to a unique letter-image, with the following differences:
\begin{itemize}
\item Here there are no identical prefixes or suffixes that need to be handled.
\item If the morphism is $q$-uniform, the distinct prefixes have length $\alpha$,
the distinct suffixes have length $\beta$, then the argument works provided that $\alpha+\beta\le q+1$.
Previously we had $\alpha=68$, $\beta=20$, and $q=95$. Here the situation is tight with $\alpha=\beta=4$ and $q=7$.
\end{itemize}
Finally, the half-flip with period $p$ in $f_3(m^\omega(\texttt{0}))$ or $f_2(m^\omega(\texttt{0}))$
implies a half-flip with period $\tfrac p7$ in $m^\omega(\texttt{0})$, which contradicts Theorem 1.1.


\begin{thebibliography}{99}

\bibitem{CuRa:2025} J.~Currie and N.~Rampersad.
Words avoiding half-flips.
\emph{International conference on combinatorics on words} (2025).

\end{thebibliography}
\end{document}